\documentclass[12pt]{article}
\usepackage{graphicx}
\begin{document}
\def\s{\underline X}
\def\y{\underline Y}
\def\S{\hat{\s}}
\def\z{{\underline Z}}
\def\ss{{\underline S}}
\def\l{\Lambda}
\def\zz{{\cal Z}^2}
\def\be{\begin{equation}}
\def\ee{\end{equation}}
\def\m{gibbs measure}
\newcommand{\ind}{\mathop{{\rm 1}\kern-0.27em{\rm I}}\nolimits}
\def\ep{\vspace{1cm}}
\parskip 2mm
\def\square{\ifmmode\sqr\else{$\sqr$}\fi}
\def\sqr{\vcenter{
         \hrule height.1mm
         \hbox{\vrule width.1mm height2.2mm\kern2.18mm\vrule width.1mm}
         \hrule height.1mm}}                  

\title{ Reconstruction of Gray-scale Images } 
\author{{ Pablo A. Ferrari \hspace{10 mm} Marco D. Gubitoso \hspace{10 mm}
 E. Jord\~ao Neves}\\ 
        Instituto de Matem\' atica e Estat\'\i stica, USP\thanks{{\it address:}
        Cx. Postal 66.281 --- 05315-970 S\~ao Paulo SP --- Brasil}}
\date{\today}
\maketitle

\begin{abstract}
  
  We present an algorithm to reconstruct gray scale images corrupted by noise.
  We use a Bayesian approach. The unknown original image is assumed to be a
  realization of a Markov random field on a finite two dimensional region $\l
  \subset \zz$.  This image is degraded by some noise, which is assumed to act
  independently in each site of $\l$ and to have the same distribution on all
  sites. For the estimator we use the mode of the posterior distribution: the
  so called {\it maximum a posteriori} (MAP) estimator. The algorithm, that
  can be used for both gray-scale and multicolor images, uses the binary
  decomposition of the intensity of each color and recovers each level of this
  decomposition using the identification of the problem of finding the two
  color MAP estimator with the min-cut max-flow problem in a binary graph,
  discovered by Greig, Porteous and Seheult (1989). {\bf Experimental
  results and a detailed  example are given in the text. We also provide a web
  page where additional information and examples can be found.}

\end{abstract}
\ep


\ep

\section{Introduction}

\label{int}

We consider here the problem of image reconstruction.  Suppose a
multicolor or gray-scale picture is subjected to noise and an observer has
access only to this corrupted version. How can she estimate the original
picture?

The analysis of this kind of problem has attracted a lot of interest
and many approaches have been considered ([G]).

One of the methods proposed is the so called MAP estimator. In this method one
assumes that the original image is a random realization of a Markov random
field that has been corrupted by some site independent noise. One assumes that
the distribution of the field (a priori distribution) is known
as well as the distribution of the noise, that is, the conditional
distribution of the observed image given the original one. The MAP estimator is
the image that has the largest probability of have produced the observed one.
This is the mode of the posterior distribution.

In multicolor images the literature proposes to use the so called Potts model
as a priori distribution. Roughly speaking, this model is a measure on the set
of images that gives more weight to images that have neighboring pixels of the
same color. This choice has a number of advantages and disadvantages. An
important disadvantage is the fact that the algorithms used to compute the MAP
estimator operate in exponential time in the number of pixels ([GJ]). 
Our main point
in this paper is to propose an alternative a priori measure with the property
that the computation time of the exact MAP estimator is polynomial in the
number o pixels. Using this approach we have produced an algorithm and a
program that reconstruct dirty images in polynomial time in the number of pixels.

In the remaining of this section we explain these ideas in some detail. In the
next section we present some {\it experimental} results obtained from the
implementation of our method to some images.  The next section is more
technical in nature. It explains why our MAP estimator can be computed
efficiently (polynomial time) and provides a proof for the theoretical result
presented below. The question of how to obtain estimators for some important
parameters associated to the observed image is discussed in the appendix.  We
close this paper with some final remarks.

In order to motivate the discussion we have to
introduce some notation. 
Assume the image is a point of $\{1,\cdots,c\}^\Lambda$, where $\Lambda$ is a
finite subset of ${\cal Z}^{2}$ (e.g. a square) with $N$ sites. The image
associates to each site or pixel in $\Lambda$ one of $c$ possible colors.
Typically values for $c$ are $c=2^8$ or $c=2^{24}$. Assume $1\le K=\log_2 c\in
\cal Z$ and let ${\underline X} = ({X}_{i})_{i = 1,...,N}$ for
${X}_{i}=(X_i^1, X_i^2, \dots , X_i^K) \in \{0,1\}^K, \;\forall\; i \in
\Lambda$, $K\ge 1$, represent the true unknown (random) image. Denote the
observed image by ${\underline Y} = ({Y}_{i})_{i = 1,...,N}$. We denote the
space of pictures $\{\{0,1\}^K\}^\Lambda$ by $\Sigma_K$.  In this notation we
have $c=2^K$ possible colors in each pixel. Each $X_i$ may correspond to a
single binary number which gives the intensity of black in pixel $i$ as
$\sum_{k=1}^K X_i^k 2^{-k}$ (in the case of gray-scale picture) or it may
correspond to three binary numbers, each giving the intensity of one of the
three {\it basic colors}. In this second case, for instance, we could have
$K=24$ and $\sum_{k=1+c8}^{(1+c)8} X^k_i 2^{-k}$, for $c = 0, 1$ and $2$
giving the intensity of, respectively, Red, Green and Blue at site $i$.  To
simplify the discussion we assume, without loss, the first case since our
approach in the second one is to reconstruct each color separately in order to
reconstruct the whole picture.

In the Bayesian setting we assume (a) that the original picture $\underline X$
is random with a known distribution which is called the {\it a priori measure}
and (b) that we know how to model the noise.

To start we consider the noise. We fix an $\epsilon >0$. 
Conditional on $\underline X$, at each bit $k$ of
each pixel $i$, independently, the observed value in bite $k$ of pixel $i$,
called ${Y}_{i}^k$, is equal to
the true value ${X}_{i}^k$ with probability $1-\epsilon$ or, with probability
$\epsilon$, it corresponds to the switched value:

\be
\label{error}
{P}_{i}({Y}_{i}^k|{X}_{i}^k) \;\; = \;\; 
\left\{ 
\begin{array}{cl} 
1-\epsilon & \mbox{if} \;\; {X}_{i}^k = {Y}_{i}^k \\  
\epsilon & \mbox{otherwise.}  
\end{array}
\right. 
\ee
Hence
\be
P({\underline Y}|{\underline X}) 
\;\; = \;\; \prod_{i \in \Lambda} \prod_{k \in \{1,\dots,K\}}
{P}_{i}({Y}_{i}^k|{X}_{i}^k) \;\; 
\ee
\be
\label{lik}
= \;\; \frac{1}{{Z}_{h}} 
\exp \left\{ h \sum_{i \in \Lambda}\sum_{k \in \{1,\dots,K\}} \mbox{\large \bf 1}_
{ \{ {X}_{i}^k = {Y}_{i}^k \}} ({\underline Y}) 
\right\}
\ee
where
\be
\label{h}
h \;\; = \;\; \log \frac{(1-\epsilon)}{\epsilon} . 
\ee
and ${Z}_{h}$ is the normalization constant.

Note that this hypothesis is not the same as the one in 
[FFG] where each pixel was either observed correctly
with probability $1-\epsilon$  or chosen uniformly among the other
 $c-1$  colors.

Now we define the MAP ({\it maximum a posteriori}) estimator.

 Denote by $\mu$ the {\it a priori} measure, that is, the distribution of
 $\underline X$. Then the MAP estimator after observing the image
  $\underline{Y}$,  $\hat{\s}= \hat{\s}(\underline{Y}) \in 
\Sigma_K$ is any
  image which maximizes the posterior distribution $P(\s|\underline{Y})
  \propto P(\underline{Y}|\s) \mu (\s)$, that is

\be
\label{map}
P(\hat{\s}|\underline{Y}) \; \; = \;\; \max_{\s \in \Sigma_K} 
P({\s}|\underline{Y})
\ee

Now we discuss the {\it a priori measure}.
We want to consider here the
situation on which  very little is known about the original
picture before the observation besides the information that it  is 
some kind of real photo, perhaps taken by a satellite, as opposed to being
a photo of something like a geometrical drawing, 
a cubist  oil  painting or a cell of a cartoon
picture. 
We suppose that the only prior information available is that the 
original picture is {\it locally smooth}. By this we mean that the measure
should be such that  neighbor pixels on the picture 
are more likely to have colors which are near in some sense. 
 More precisely we assume that there exists a 
real valued function $H$ which  
indicates how {\it
smooth } is an image and  take the {\em a priori
measure} to be:

\be
\label{gibbs}
\mu(\s) = \frac{e^{-\beta H(\s)}}{Z}
\ee 
where $Z = \sum_{\s}{e^{-\beta H(\s)}}$, with the sum taken over all possible
images,  is a normalization constant and
${\beta} \geq 0$ is a real parameter.  This parameter measures 
the tendency to be smooth since if $\beta$ is large $\mu$ is concentrated
on images with small values of $H$ while if $\beta$ is small $\mu$ is close to
the uniform distribution on $\Sigma_K$, the set of all images.
 The motivation for this formula, 
in particular the minus sign in front of $\beta$,
 comes from statistical mechanics where it would be called Gibbs
measure, $\beta$  would be the inverse temperature and $H$ would be
the Hamiltonian (or Energy) function.

Plugging (\ref{lik}) into (\ref{map}), and taking logarithms, we get that the
images which maximize the posterior distribution  (\ref{map}) are those which
maximize
\be
\label{loglike}
\beta H(\underline X) + h \sum_{i \in \Lambda}\sum_{k \in \{1,\dots,K\}} \mbox{\large \bf 1}_
{ \{ {X}_{i}^k = {Y}_{i}^k \}} ({\underline Y}) 
\ee 

The image that maximizes this expression 
makes the best compromise between being globally {\it smooth},
regulated by $\beta H $,  and
agreeing as much as possible with the observed image which is regulated by
$ h \sum_{i \in \Lambda}\sum_{k \in \{1,\dots,K\}} \mbox{\large \bf 1}_
{ \{ {X}_{i}^k = {Y}_{i}^k \}} ({\underline Y}) $. 
Of course there exists only one relevant parameter in the
maximization problem, say $h/\beta$.

Although {\it any} strictly positive measure $\mu$ on $\Sigma_K$ can be
written as above {\it for some } $H$, and thus (\ref{gibbs}) can be safely
assumed for a generic {\it a priori } measure, we are interested only in the
case on which $H$ is both {\it simple} and {\it reasonable } as a measurement
of smoothness. This simplicity requirement will basically mean the
assumption
that $H$ has only local dependence as follows

\be
\label{genham}
H(\s) = \sum_{<i,j>} d(X_i, X_j)
\ee
where $d(X_i, X_j)$ is a measure for 
 the {\it distance} between the colors at pixel
$i$ and $j$ and 
 the sum is taken over all pairs of {\it  neighbor} pixels  in the
 picture. That is, $<i,j>= \{(i,j): \vert i-j\vert =1\}$.

    The  choice of $d(X_i, X_j)$  depends on what kind of local
properties one expects in the original picture.

One choice which is common in the literature  (see, for instance [FFG] and its
quotations) is
\be
\label{potts}
H_P(\s) = \sum_{<i,j>} \ind_{\{X_i \neq X_j\}}
\ee
where the sum is taken over all pairs of neighbor sites  in the lattice. 

This choice is related to the Potts model in statistical mechanics (see for
instance [M]) and explains the subscript. It is a model with very interesting
properties that  corresponds, in the
two-color case, to the Ising model [MW]. To find a solution to
(\ref{map}) we have to find what is called in the physics literature a {\it
  ground state} for the Ising model with random magnetic field. This field is
induced by the observed image.

 Most of the interesting properties of those statistical mechanics
 models, like {\it phase transition}, 
appear in the so called {\it thermodynamic limit}, the limit on which the size
of the system grows to the whole lattice ([R], [MW]).
Even though in our discussion the lattice
size is kept fixed we will need in section 4 some exact results on the thermodynamic
 limit in the two-color case (Ising model).
Also some finite size considerations like the effect of boundary
conditions (in our case we chose {\it free} boundary conditions) may be
important.

The measure defined by  (\ref{gibbs}), with  (\ref{potts}) plugged in,
 would 
 describe  images on which neighbor pixels tend to be equal but
if two pixels have different colors 
the cost of this {\it interface}  does not
depend on {\it how different } they are.
If two
pixels have different colors the most likely is that each one belongs
to a single-color region  separated by a {\it sharp} interface.

We could also try to represent the situation on which this is not
necessarily the case choosing  a  
finer notion of 
 distance between $X_i$ and
$X_j$. Two somewhat natural choices would be

\be
\label{dist}
H_1 (\s) = - \sum_{<i,j>} \left|\sum_{k=1}^{K} (X_i^k - X_j^k)\frac{1}{2^k}\right|
\ee
 
or 

\be
\label{dist2}
H_2 (\s) = - \sum_{<i,j>} \left(\sum_{k=1}^{K} (X_i^k -
  X_j^k)\frac{1}{2^k}\right)^2. 
\ee

Note that, for any choice of  $H$,
 the MAP problem is well posed since we could, at
least in principle, check the finitely many values in the right hand side of 
(\ref{loglike}) and choose a image that maximizes it. But 
since $|\Sigma_K| = {2^K}^N$
with something like $K=24$ and $N= 400 \times 600 $ this approach is not
computationally feasible. The problem with all the above  choices 
 is that it is not known how to do better
than this time-consuming 
maximization by inspection in the multicolor case ($K>1$) and thus
the problem is, in practice, not solvable.

We propose  another choice for $H$, intermediate
between the Potts and the other two mentioned above, which 
respects nicely the notion of
smoothness which led to the choice (\ref{dist}) but nevertheless  induces
 to a polynomial time
maximization problem.  This choice is:

\be
\label{H}
H (\s) = \sum_{k=1}^K  \sum_{<i,j>} |X_i^k - X_j^k|\frac{1}{2^k}=\\
 \sum_{k=1}^K \sum_{<i,j>}  \ind_{\{X_i^k \neq X_j^k\}}\frac{1}{2^k}
\ee

With this, the problem of finding the $2^K$-color 
 image which maximizes the 
posterior distribution $P(\s|\underline{Y})
  \propto P(\underline{Y}|\s) \mu (\s)$ is decomposed into 
 $K$ binary color
maximization problems. Namely one has to solve 

\be
\label{mapk}
P(\hat{\s}_k|\underline{Y}_k) \; \; = \;\; \max_{\s \in \Sigma_2} 
P(\underline{Y_k}|\s) \mu_k (\s)
\ee
for each $k$,  $1 \leq k \leq K$ where  $\underline{Y}_k= \{Y_i^k\}_{\{i\in
  \Lambda\}}$ 
is the $k$-th component of the observed image and $\mu_k$ is the Gibbs measure
defined on the space of binary color images, $\Sigma_2$, by (\ref{gibbs}) with 

\be
\label{hk}
H^{k}(\underline{X}_k) = \sum_{<i,j>}  \ind_{\{X_i^k \neq X_j^k\}}\frac{1}{2^k}. 
\ee

The $K$ color image given
by $\underline{Y}$ is a solution of (\ref{map}). Each one of these $K$ binary
images is called a {\it layer}. The weighted sum of the MAP for each layer gives our
(gray-scale) image estimator:
$$
\hat X_i = \sum_{k=1}^K 2^{-k} \hat X^k_i .
$$

On the multicolor case one solves a  {\it gray-scale}
problem for each basic color.

Each binary problem can be solved in polynomial time using the results by
Greig, Porteous and Seheult who reformulated it 
as one involving finding a minimum cut on a capacitated network [FF]
for which there exist fast algorithms. These ideas are presented 
briefly in the next section.
 From the statistical point of view the approach is different whether
the parameters $\beta$ and $\epsilon$ are known or not. A truly Bayesian
approach would associate an a priori measure to each one of the
parameters $\beta$ and $\epsilon$. We leave this alternative to future
work.  Another possibility is to estimate these parameters from the
observed image using classical frequentist analysis.  This is possible
under the hypothesis on the noise and on the original image being a
sample of a Gibbs measure for the Potts Hamiltonian (\ref{potts}) (which
is always the case after the decomposition in binary colors) using exact
results for the two-dimensional Ising model obtained by Frigessi and
Piccioni [FP]. We have produced an algorithm based on [FP] for
estimating the parameters $\epsilon$ and $\beta$. This is explained in
section 4 below.

Once one has a fast algorithm to reconstruct images it is natural to ask what
is the effect of iterating the whole process. More precisely, what happens if
one takes the reconstructed image and applies the method again, perhaps
updating the value of $h/ \beta$?
  
Note that to apply the method again for the (already) reconstructed image is
equivalent to assume that the reconstructed image could be thought as obtained
 from some original picture which was chosen with respect to $\mu$ and then
subjected to noise.  Even though it is not difficult to verify that this
assumption is false the question is interesting both from the mathematical and
from the applied point of view. On the mathematical side one has a mapping,
$\underline Y \mapsto \hat{\s}(\underline Y)$, from $\Sigma_K$ into itself, a
discrete time dynamical system, and it is natural to ask about iteration
properties.  On the practical minded side one can ask about the effect of
iteration on the {\it quality} of reconstruction even if this can only be judged
subjectively.

If we decide to update $h/\beta$ before each iteration the natural thing to do
would be to chose a larger value at each step since we assume the procedure
did a good job and removed some of the noise, therefore $\hat{\s}(\underline
Y)$ would be more reliable than $\underline Y$. It would then be natural to
increase $h$, without changing $\beta$. A natural procedure to find each
updated value is to use again the estimators given by Frigessi and Piccioni
([FP]). In doing this we find, in fact, that this parameter increases, after
one iteration, as expected.

What is rather surprising is that once one tries this iteration procedure,
either keeping $h/\beta$ fixed or increasing it, at each step, one finds
that {\it it has no effect et all}. The twice reconstructed image is {\it exactly}
equal to the once reconstructed one. In other words, $\hat{\s}$ as a function from
$\Sigma_K$ into itself has a fixed point in each once-reconstructed image.

More precisely write $F_{(h/\beta)} ({\s})= \hat{\s}$ for
 the function from $\Sigma_{K}$
into itself defined in (\ref{map}). Then we have the following

\noindent {\bf Proposition:} $F_s({\s}) = F_t (F_s(\s))$ for all $t \geq s$.
 
We prove this proposition in the next section after reviewing 
some
results on networks and the connection with the maximization
problem considered in this note.

Given this result a natural question arises. Is this
property not true {\it in general}? Suppose one has a random function
assuming values in some space $S$ with some unknown parameter ${\it s}$ which
itself belongs to $S$ chosen according to some {\it a priori} measure 
 (in our case, $\it s$ corresponds to the original
image). Define the MAP estimator as usual to obtain a function from $S$ into
itself. Does it always satisfy the this fixed point property? 
One could argue that since all the 
information about the original image which is  contained in the observed one
should be still present in the reconstructed image then this iteration should give
no further results and the {\it twice} reconstructed image should  {\it
  always} be equal to the  
once-reconstructed one. As it turns out this is
not the case and it is not difficult to find counterexamples.


\section{Experimental Results}

We developed a program (\texttt{map}) in \emph{C} to reconstruct images
according to the method proposed. This program reads an image in the
\emph{portable bit map} format for true-color (\emph{.ppm}) or gray-scale
(\emph{.pgm}) pictures.  The intensity of each color in each pixel (or the
single gray scale intensity there) is written as a (eight bits) binary number.
The program then constructs a graph for each one of those layers (as defined
after equation (\ref{hk})), construct the corresponding graphs and finds the
minimum cut for each one of them using a standard implementation of the Ford
Fulkerson algorithm. The program uses integer arithmetic scaled by a factor of
10000.  The values for $\beta$ and $H$ are supplied by command line switches.
It is also possible to supply a bit mask to select which planes to
reconstruct.  The program is portable and it has been tested on \emph{SunOS},
\emph{Solaris} and \emph{Linux}.

The main problem in developing the program was the large amount of
memory to hold the image, the graph and the auxiliary data structures.
The solution was to use a compact representation of the graph, namely a
matrix where each cell has a pixel and the values for the flow in each
of 5 directions (4 neighbors plus $s$ or $t$). This avoided the explicit
representation of the edges and the extra space for the image itself.

A series of tests were made on a 700Mhz \emph{Pentium III} machine, with
128Mb of main memory, running Linux Debian potato 2.2, kernel version
2.4.  The running times for restoration of six sample images are
presented in table \ref{times}.  These images are shown in tables
\ref{img1}, \ref{img2} and \ref{img3}.

The images were obtained from a picture shot with a digital camera and
converted to 8 bit gray-scale.  The final result was saved in the file
``\texttt{fish.pgm}''.  Two aditional files were obtained by clipping
the image to a rectangle with half of the area (``\texttt{fish2.pgm}'')
and to another one with a quarter of the total area
(``\texttt{fish3.pgm}'').

The noise was added ``artificially'' by a program which scan every bit
in the image and inverts it with a given probability. The probabilies
used were 10\% and 15\%.  The modified files were named
\texttt{fish.10.pgm}, \texttt{fish2.10.pgm}, \texttt{fish3.10.pgm} and
\texttt{fish.15.pgm}, \texttt{fish2.15.pgm}, \texttt{fish3.15.pgm},
respectively.

Each image was restored in all bit planes using 3 values of $\beta$,
namely 0.1, 0.3, 0.5. An aditional restoration, just on the most
significant bit plane were made using the estimated $\hat\beta$ as
indicated in the appendix.  The values are presented in the tables
\ref{img1}, \ref{img2} and \ref{img3}.

We did not try to define a metric to measure the quality of the
reconstruction, relying on a subjective analysis.  Besides the removal
of noise, it was observed an improvement in the shades and
(consequently) in the third dimension perception, even if some blurring
is introduced.

These and other examples are available directly from the authors or at
  \verb+http://www.ime.usp.br/~gubi/MAP/index.html+
\nopagebreak[4]
\begin{table}[h!]
  \begin{center}
    \begin{tabular}[c]{|l|rrr|}
      \hline
      & \multicolumn{3}{c|}{Beta}\\
      \raisebox{1.5ex}[0pt]{\large File} & 0.10  &0.30  &0.50\\\hline
      \texttt{fish.10.pgm}  &6.40  &17.77  &16.78  \\
      \texttt{fish.15.pgm}  &6.05  &14.38  &19.54  \\
      \texttt{fish2.10.pgm} &2.81  & 4.45  & 6.35  \\
      \texttt{fish2.15.pgm} &2.65  & 4.38  & 5.86  \\
      \texttt{fish3.10.pgm} &1.38  & 2.43  & 3.20  \\
      \texttt{fish3.15.pgm} &1.34  & 2.37  & 3.02  \\
      \hline
    \end{tabular}
  \caption{Restoration running times for three sample images with two
    noise levels (10\% and 15\% per bit)}
  \label{times}
  \end{center}
\end{table}

\begin{table}[h!]
  \begin{center}
    \begin{tabular}[c]{|c|c|c|}
      \hline
      Original & 10\% & 15\% \\
      \hline
      \includegraphics{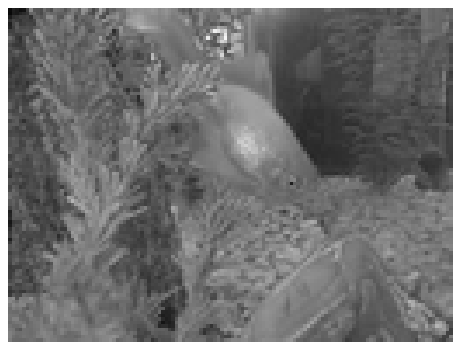} &
      \includegraphics{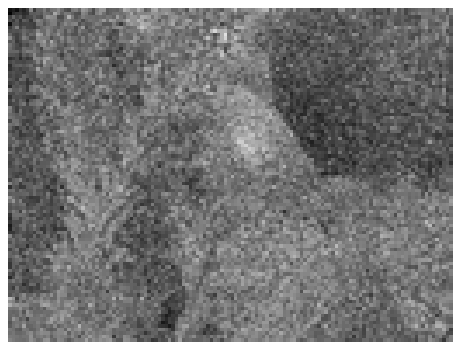} &
      \includegraphics{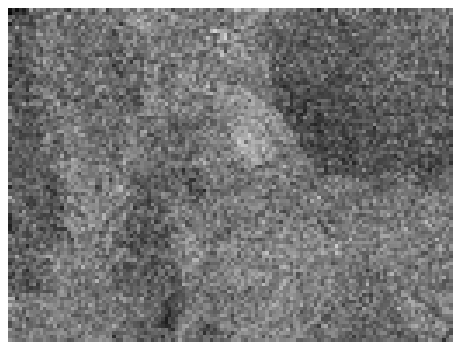} \\
      \hline
      \raisebox{1.6cm}{$\beta=0.1$}&
      \includegraphics{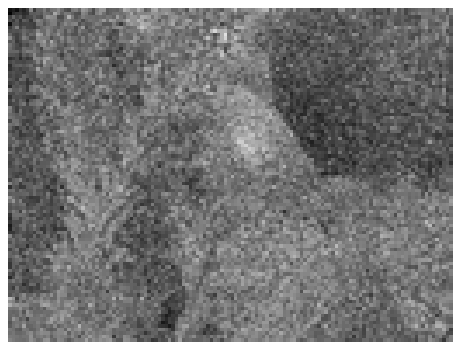} &
      \includegraphics{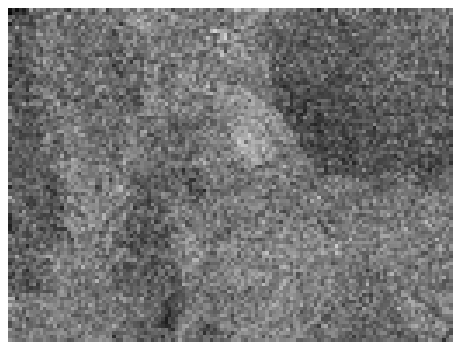} \\
      \hline
      \raisebox{1.5cm}{$\beta=0.3$}&
      \includegraphics{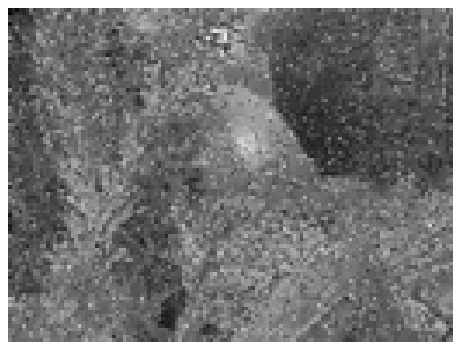} &
      \includegraphics{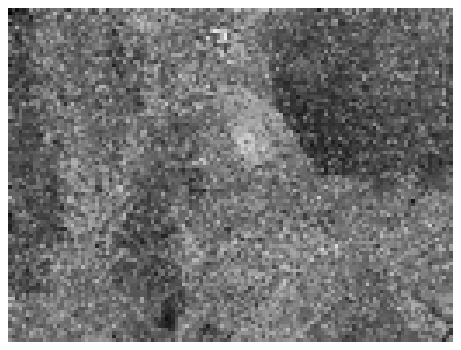} \\
      \hline
      \raisebox{1.5cm}{$\beta=0.5$}&
      \includegraphics{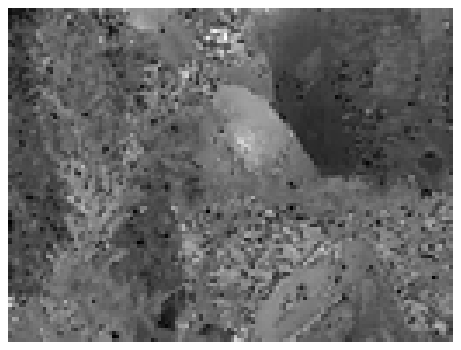} &
      \includegraphics{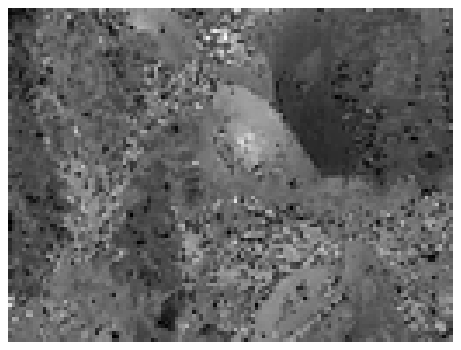} \\
      \hline
      \raisebox{1.5cm}{$\hat\beta_{10\%}=\hat\beta_{10\%}=0.45$}&
      \includegraphics{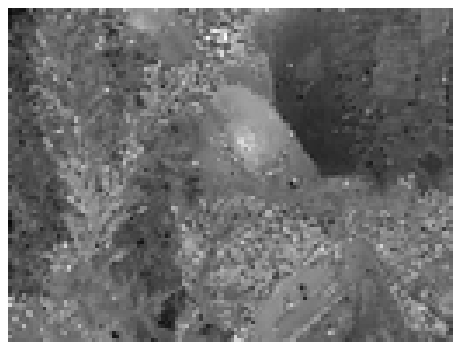} &
      \includegraphics{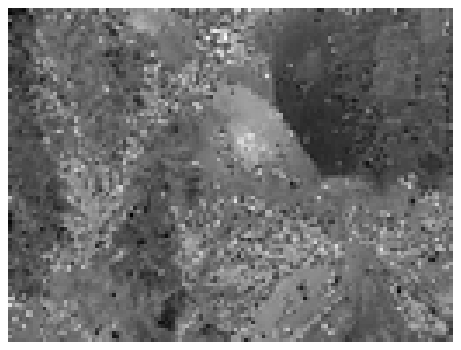} \\
      \hline
    \end{tabular}
    \caption{Sample image (\texttt{fish.pbm})}
    \label{img1}
  \end{center}
\end{table}
\clearpage
\begin{table}[h!]
  \begin{center}
    \begin{tabular}[c]{|c|c|c|}
      \hline
      Original & 10\% & 15\% \\
      \includegraphics{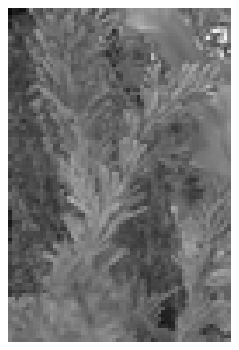} &
      \includegraphics{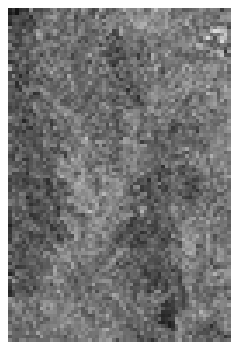} &
      \includegraphics{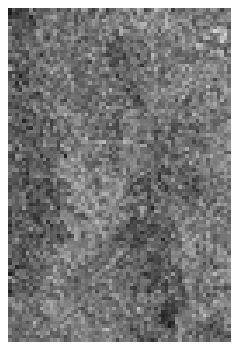} \\
      \hline
      \raisebox{1.5cm}{$\beta=0.1$}&
      \includegraphics{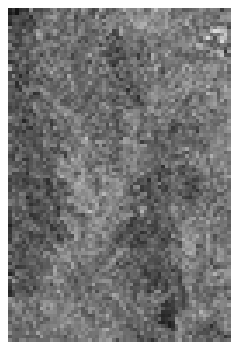} &
      \includegraphics{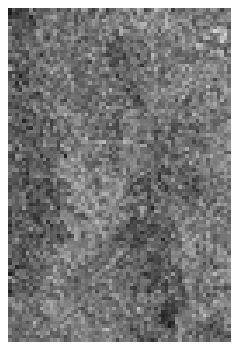} \\
      \hline
      \raisebox{1.5cm}{$\beta=0.3$}&
      \includegraphics{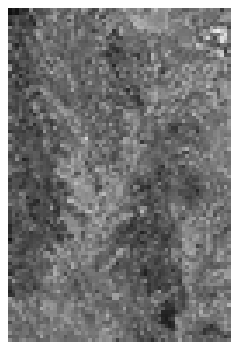} &
      \includegraphics{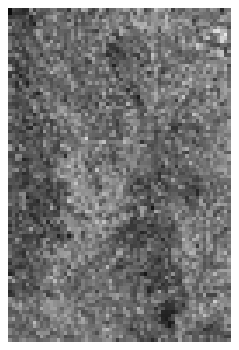} \\
      \hline
      \raisebox{1.5cm}{$\beta=0.5$}&
      \includegraphics{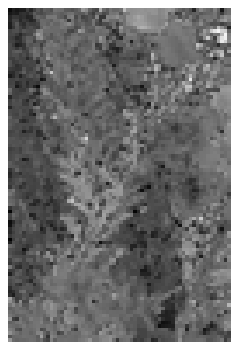} &
      \includegraphics{fmeio.10.limpo.5.ps} \\
      \hline
      \raisebox{1.5cm}{$\hat\beta_{10\%}=\hat\beta_{15\%}=0.45$}&
      \includegraphics{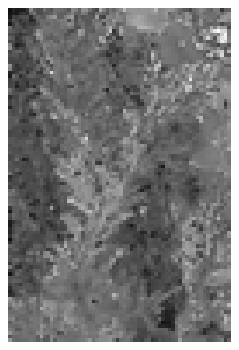} &
      \includegraphics{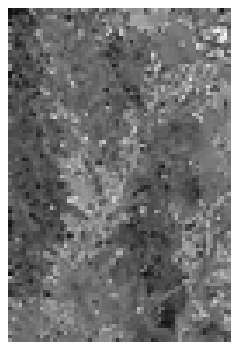} \\
      \hline
    \end{tabular}
    \caption{Same photo of table \protect{\ref{img1}}, but cut in half
      (\texttt{fish2.pbm})}
    \label{img2}
  \end{center} 
\end{table}
\clearpage
\begin{table}[h!]
  \begin{center}
    \begin{tabular}[c]{|c|c|c|}
      \hline
      Original & 10\% & 15\% \\
      \includegraphics{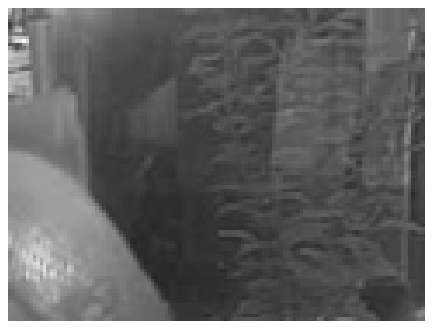}&
      \includegraphics{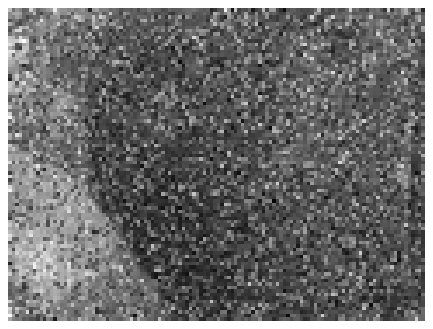} &
      \includegraphics{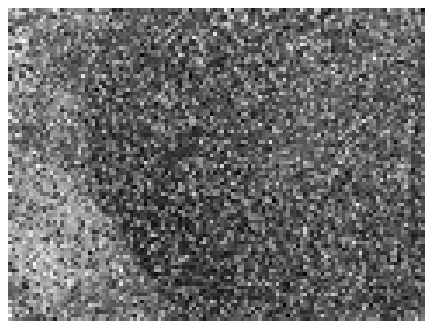} \\
      \hline
      \raisebox{1.5cm}{$\beta=0.1$ }&
      \includegraphics{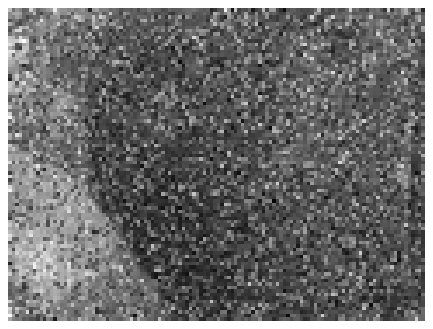} &
      \includegraphics{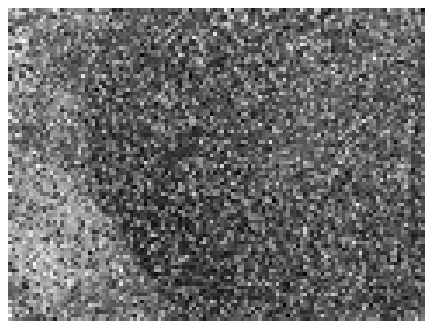} \\
      \hline
      \raisebox{1.5cm}{$\beta=0.3$}&
      \includegraphics{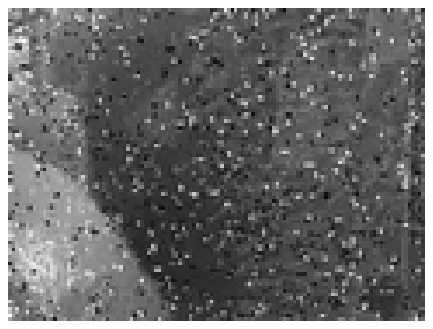} &
      \includegraphics{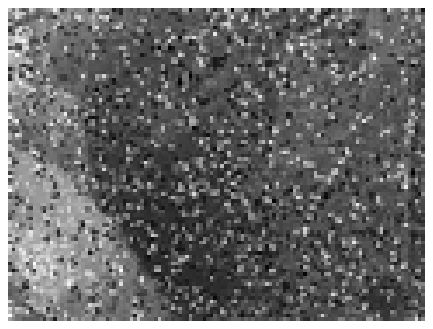} \\
      \hline
      \raisebox{1.5cm}{$\beta=0.5$}&
      \includegraphics{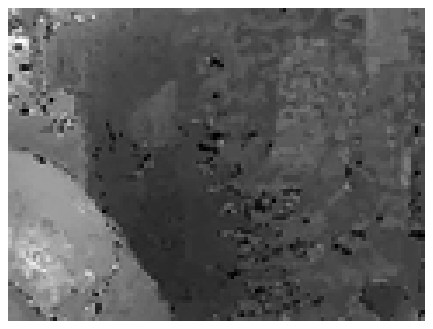} &
      \includegraphics{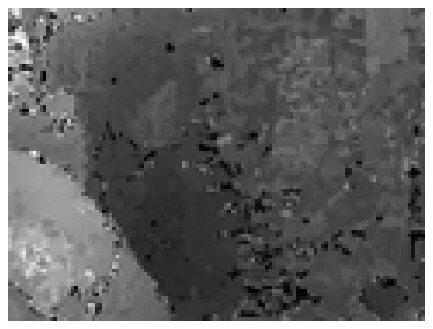} \\
      \hline
      \raisebox{1.5cm}{$\hat\beta_{10\%}=0.549\quad
        \hat\beta_{15\%}(0.590)$}&
      \includegraphics{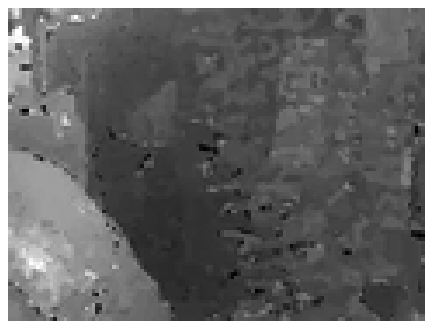} &
      \includegraphics{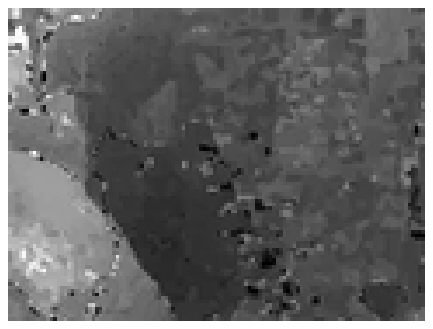} \\
      \hline
    \end{tabular}
    \caption{$\frac{1}{4}^{\mbox{th}}$ of same image (\texttt{fish3.jpg})}
    \label{img3}
  \end{center}
\end{table}
\clearpage
  

\section{MAP estimator and networks}

In this section we present some capacitated network ideas
which are used to solve the two-color maximization problem described before as
discovered  by Greg, Porteous and Seheult [GPS] 
and use them to prove Proposition. 


As mentioned before with the choice 
of $H$ given by (\ref{H}) the maximization in (\ref{map}) is
decomposed into $K$ two-color problems. 
Therefore we assume in this section
that we are in the binary color case $K=1$.

The MAP estimator given $\y$  is the image $\hat{\s} \in \{0,1\}^{\Lambda}$
 which gives the maximum of 
\be
\label{L}
L_{\alpha}({\underline{X}}|{\underline{Y}})=
 \alpha \sum_{i \in \Lambda} \mbox{\large \bf 1}_
{ \{ {X}_{i} = {Y}_{i} \}} ({\underline X}) 
+  \sum_{<i,j>}  \mbox{\large \bf 1}_{\{X_i = X_j\}} ({\underline{X}})
\ee
where $\alpha = h / \beta$.

A network $N$ is a graph $G= (V, E)$, where $V$ is a finite set of vertices
and $E$ is a set of couples of vertices, with a capacity $c(e)\geq 0$
associated with each edge $e \in E$.

We define networks on the set of vertices $V$,
given by the sites in $\Lambda$ plus two extra ones
denoted by $s$ (source) and $t$ (sink), that is
$$
V = \Lambda \cup \{s\} \cup\{t\}.
$$
 For the set of arcs we take 
$$
E =  \left(\cup_{i}\{(s,i)\}\right)\cup\left(\cup_{i} \{(i,t)\}\right)
\cup\left(\cup_{<i,j>}\{(i,j)\}\right)
$$ 
where the first two unions are taken on $i\in \Lambda$ and the last one
is taken over pairs of nearest neighbor sites in $\Lambda$. 

Given the observed image $\underline{Y}$ and a real number $\alpha$ we define
the capacities of the network $N_{\alpha} (\y)$ as follows.

If $Y(i)=1$ we set $c_{\y,\alpha}(s,i)=\alpha$ as its capacity,  
 otherwise set 
$c_{\y,\alpha}(i,t)=\alpha$; to each arc $e=(i,j)$ of neighbor sites in 
$\Lambda$ we associate $c_{\y,\alpha}(e)=1$. All other arcs have capacity
 zero.

For each  image $\underline{X}$ let 
$$A(\underline{X})= \{ s\} \cup \{ i \in
\Lambda: 
X(i)=1\}$$ 
$$B(\underline{X}) =\{t\}\cup \{  i \in \Lambda : X(i)=0\}.$$ 
 These two sets define
a {\it cut} of the network 
$$
{\cal C}(\s)=  \{ (i,j) \in E: i \in A(\s), j \in B(\s)\}. 
$$
Notice that the cut $ {\cal C}(\s)$ consists of a set of arcs whose removal
(cut) makes it impossible to find a path going from $s$ to $t$ through 
arcs with non-zero capacity.  
If we now define the {\it capacity of the cut} $ {\cal C}(\s)$ by the quantity
\be
\label{cut}
C_{\y,\alpha}(\underline{X})= 
\sum_{(i,j) \in {\cal C}(\s)} c_{\y,\alpha}(i,j)
\ee
it is very simple to check that 
$$
L_{\alpha}(\underline{X}|\underline{Y}) = a - C_{\y,\alpha}(\underline{X}),
$$ 
where $a$ is a constant which
does not depend on $\underline{X}$. Therefore to find the MAP estimator one
has to find the cut which minimizes (\ref{cut}): the so called
{\it minimum cut} [FF].
 
Ford and Fulkerson showed that the value of the capacity of the minimum cut is
equal to the {\it maximum flow } through the network from source to sink.
Recall that a flow ${\bf f}$ in a network $N$ (on $G(V,E)$ with capacity
$\{c(e)\}_{e\in E}$) from $s$ to $t$ is a collection of real numbers $\{ {\bf
  f}(e)\}_{e \in E}$, where ${\bf f}(e)$ can be thought as the amount of fluid
per unit time going through the {\it pipeline } $e$, such that the flow on
each arch does not exceed its capacity, $0\leq {\bf f}(e) \leq c(e)$, for all
$e \in E$, and such that the flow is conserved

\be
\sum_{j \in V}{\bf f}(i,j) - \sum_{j \in V} {\bf f}(j,i) = 0
\ee
for all $i \ne s$.

Well known fast algorithms to find this maximum flow exist and therefore the
problem is solved in practice. The algorithm presented in this paper uses this
method to reconstruct each one of the (8 or 24) layers.

\vskip 3mm

{\noindent\bf Proof of the Proposition.}
To simplify the notation write $\hat{\s}=F_{\alpha}({\underline{Y}})$.
Let $\alpha' \geq \alpha$. 
The networks defined by the pair 
$(\underline{Y}, \alpha)$, denoted by  $N_{\alpha} (\y)$,  and by 
$(\hat{\s}, \alpha')$, $N_{\alpha'} (\S)$,
can   differ only in the capacities assigned to each arc.

The proposition  asserts that the cut defined by $\S$ also has
minimum capacity in $N_{\alpha'} (\S)$, that is we want to prove that
\be
\label{*}
C_{\hat{\s},\alpha'}(\hat{\s}) \leq 
C_{\hat{\s},\alpha'}(\underline{Z}),
\ee
for any image $\underline{Z} \in \Sigma_1$.

Fix  $\underline{Z}$ and let $E_I ={\cal C}(\S) \setminus  {\cal C}(\z)$, 
$E_{II}={\cal C}(\S) \cap  {\cal C}(\z)$
and $E_{III}=   {\cal C}(\z) \setminus  {\cal C}(\S)$.
By definition of MAP estimator  we have
\be
\label{inalpha}
C_{\y,\alpha}(\hat{\s}) \leq 
C_{\y,\alpha}(\underline{Z}),
\ee
which implies
\be
\label{1}
D_{\y,\alpha}(E_{I}) \leq 
D_{\y,\alpha}(E_{III}).
\ee
where, for a set of edges $E$, 
$$
D_{\y,\alpha}(E) = \sum_{e\in E} c_{\y,\alpha}(e)
$$

Therefore to check (\ref{*}) it is enough to verify
\be
\label{2}
D_{\S,\alpha'}(E_{I}) \leq D_{\y,\alpha}(E_{I})
\ee
\be
\label{3}
D_{\S,\alpha'}(E_{III}))  \geq D_{\y,\alpha}(E_{III}).
\ee

We start with inequality (\ref{2}). Assume $e \in E_{I}$. If $e$ is an
external edge (connecting some site $i$ with either the source $s$ or the sink
$t$), then, since $e\in {\cal C}(\S)$, it must be $c_{\S,\alpha'}(e)=0$ which
is less or equal than $c_{\y,\alpha}(e)$. 
On the other hand, if $e$ is an internal edge (i.e. $e = (i,j)$ for $i$ and $j$
nearest neighbors in $\Lambda$)  then its capacity
equals $1$ in both cases. This concludes the proof of (\ref{2}).

Suppose now $e \in E_{III}$. Inequality (\ref{3}) follows  from the observation
that $c_{\S,\alpha'}(e)$ can not be zero. \square


\section{Final Remarks}

The different approaches used in image reconstruction are based in quite
different set of theoretical ideas and it is not clear how to compare their
results.  One possible measure for the {\it quality } of the reconstruction,
used in [FP] to compare 9 algorithms, is to evaluate the proportion of
pixels classified correctly. Since our main goal here was to present a method
for the reconstruction of multicolor images we leave the comparison with other
methods for future work. In any case, since we are working with exact MAP's
for the chosen Hamiltonian, our method will be as good (and as bad) as the
usual two colors MAP estimators, regarding the proportion of bits reconstructed.

For the usual MAP reconstruction problem in the multicolor case no fast
algorithm is known ([FFG]). For a probabilistic approach via {\it Simulated
  Annealing} in order to get the exact estimator one needs to decrease very
slowly some parameter while the computation goes on and thus needs a
prohibitively large amount of time [Gi].

One possibility is to accept approximated MAP estimators which can be obtained
fast enough. One can do this with simulated annealing by updating the
parameter fast enough but then we lose control on how close the approximation
is to the exact one.  An approximate of the MAP estimator with a probabilistic
analysis of the error in the three color case was developed in [FFG].
Another approach can be found in [J].

Our approach is not completely Bayesian as we also consider a situation on
which the {\it a priori} measure has some unknown parameters but do not assume
any prior knowledge about them. In the Appendix we describe a method to
estimate these parameters from the picture itself using {\it classical
  statistical methods}, as proposed by Frigessi and Piccioni [FP].  In those
cases we verify that plugging those estimated values into the formulas used to
get the MAP estimator provides a reconstruction which appears to be the best.
As mentioned before we do not try to quantify this.

Ricardo Maronna proposed that instead of estimating $\beta$ and $h$ one could
look for the $\alpha = \beta/h$ that maximizes (\ref{loglike}). The estimator
for the true image will then be the $\s$ which realizes this maximum with the
best $\alpha$. It is not clear yet how to justify this theoretically. Another
alternative would be to choose the uniform distribution on an appropriate
range as a priori distribution for $\epsilon$ and $\beta$. However the
computation of the MAP estimator in this case seems prohibitive.

As an experimental observation we remark that reconstructing only the first layer
of a dirty image (and leaving the others as they are) gives a quite good
visual result. A possible explanation of this fact is that, as a consequence
of the binary decomposition, each layer is ``half'' as important as the
previous one.

Looking for algorithms which give smooth solutions we propose the following
hierarchical procedure. First consider layer $1$ and find $\hat Z^1$, as
before for $\hat X^1$. Then given layer $\hat Z^{\ell}$, for
$\ell=1,\dots,k-1$, define $\hat Z^{k}$ as the (binary) image that maximizes
\be
\label{LL}
L_{\alpha}({\underline{X}_k}|{\underline{Y}})=
 \alpha \sum_{i \in \Lambda} 
\mbox{\large \bf 1}_
{ \{ {X}^k_{i} = {Y}^k_{i} \}} ({\underline X}) 
+  \sum_{<i,j>}  \prod_{\ell =1}^{k} \mbox{\large \bf 1}_
{ \{ {X}^\ell_{i} = {X}^\ell_{i} \}} ({\underline{X}}).
\ee
In words, this algorithms will try to get the same value for two neighbors in
layer $k$ if these neighbors have the same value for all previous layers. If
not they are not coupled. Visual realizations of this algorithm give also good
results. Although there is no Hamiltonian for this model, the algorithm is
 well defined. Moreover each layer corresponds to the so called diluted Ising model.

Assume that we get $N$ independent samples of the same image $\underline
Y_1,\dots,\underline Y_N$. 
This means that the
original realization $\s$ of the image is the same but the noises are
independent. A generalization of our approach deals with this problem
in the following way. We modify the capacities associated to the graph
$G$. To each arc $(s,i)$ connected to the source we associate the capacity 
$$
\alpha \left(\sum_{n=1}^N (2 Y^k_{n,i}-1)\right)^+.
$$
and then to each arc $(i,t)$ connected to the sink we associate the capacity 
$$
\alpha \left(\sum_{n=1}^N (2 Y^k_{n,i}-1)\right)^-.
$$
In the same vein one could use the remaining two colors to get information
about the color being reconstructed.

\parskip 2mm
\parskip 2mm


\section{{\bf Appendix: } Estimation of $\beta$ and $\epsilon$}

In this section we apply  ideas from Frigessi and Piccioni ([FP])
which exploit well known (but highly nontrivial) results from
statistical mechanics ([R], [MW]) to
obtain estimators for 
the parameters $\beta$ and $\epsilon = (1+e^h)^{-1}$.

The result is a
program called {\bf estima} that has as input a multicolor image and
returns one estimated $\beta$ and one estimated $\epsilon$ for each
layer.

The measure defined in (\ref{gibbs}) with $H$ given by
(\ref{potts}) in the two-color case (Ising model) in the
thermodynamic limit gives rise to a translation invariant measure
which may be ergodic or not according to the value of $\beta$.

This choice of $H$ does not favor zeros or ones. This symmetry is
 perhaps easier to see if we represent a configuration as an element  
$\underline{S}$ in
 $\{-1, +1\}^{\zz}$ with $-1$ replacing zeros. This is  the usual Ising
 notation while the one using zeros and ones is known in statistical
 mechanics as {\it lattice gas} representation. In the situation
 considered here they are equivalent. 

More precisely in if 
 $\underline{Y} \in \{0,1\}^{\Lambda}$ let $\ss$ be the corresponding
 Ising configuration given by

\be
\ss= \{{S}_{i}\,,\, i \in {\Lambda}: S_i = 2 Y_i -1\}
\ee

And instead of Hamiltonian (\ref{potts}) we use

\be
\label{ising}
H_P(\ss) = - \sum_{<i,j>} S_i  S_j
\ee

The corresponding infinite volume limit measure is 
not ergodic if $\beta$ is smaller than some (known) critical
value $\beta_{c}$. In this case the limit measure is a mixture
with equal weights of 
two measures $\mu^{+}_{\beta}$ and $\mu^{-}_{\beta}$, the first favoring 
configurations with more $+1$'s  than $-1$'s and the other favoring 
configurations with more $-1$'s than $+1$'s. If $\beta \geq \beta_c$ then
$\mu^{+}_{\beta} =\mu^{-}_{\beta}$. If we denote by 
$E_{\beta}^+ (f)$ ( $E_{\beta}^- (f)$ ) the expected value of a
function $f$ defined on  $\{-1, +1\}^{\zz}$ the symmetry 
between $\mu^{+}_{\beta}$ and $\mu^{-}_{\beta}$ imply

\be
E_{\beta}^+(S_i S_j)=E_{\beta}^-(S_i S_j)= r_{\beta}(|i-j|)
\ee

 $r_{\beta}(|i-j|)$ is called the two point correlation function for
the infinite volume system at inverse temperature $\beta$. We will
need two of those correlation functions in what follows:
$r_{\beta}(1)$, for nearest neighbors, and $r_{\beta}(\sqrt{2})$ for
neighbors along the diagonal on the lattice.

 Suppose $\ss$ corresponds  (in the
Ising notation) to the
original image and $\underline{R}$ corresponds to the observed image after
noise. Under the hypothesis on the noise

\be
E_{\beta}^+(R_i R_j)=E_{\beta}^+(S_i S_j) (1-2\epsilon)^2
\ee

Therefore the ratio between $r_{\beta}(1)$ and $r_{\beta}(\sqrt{2})$
computed for the observed image depends only on $\beta$. Call this
ratio $\phi(\beta)$.


 
If $\underline{R} \in \{-1,1\}^{\Lambda}$
corresponds to an observed  two-color picture [FP] found a sequence
$(\hat{\beta},\hat{\epsilon})$ of consistent estimators for $(\beta,
\epsilon)$ given by 

\be
\hat{\beta} = {\phi}^{-1} \left(\frac{G^1(\underline{R})}{G^2(\underline{R})}\right)
\ee
and 
\be
\hat{\epsilon}=\frac{1}{2}\left\{1-\left(\frac{G^1(\underline{R})}
{r_{\hat{\beta}}(1))}\right)^{\frac{1}{2}}\right\}
\ee
with 
\be
G^a(\underline{R})= \frac{\sum^a_{(i,j) \in \Lambda^0} R_i R_j}{4|\Lambda^0|} 
\ee
where for a =1 the sum is taken over pairs of nearest neighbor sites along the
 lattice directions in $\Lambda^0$, the interior of $\Lambda$, and for a =
2 the sum is over neighbor sites along the two lattice diagonals again in the
interior of $\Lambda$.

Expressions for  $r_{\beta}(1)$ and $r_{\beta}(\sqrt{2})$ are 
complicated involving elliptic integrals with  different formulas for
$\beta$ smaller, equal and larger than $\beta_c$ and the function 
 ${\phi}^{-1}$ must be
computed numerically. 
Our program {\bf estima} finds these estimators from a observed image.

\noindent{\bf Acknowledgments.}

We thank Ricardo Maronna for comments on a draft of this paper.
We thank Eric Soares da Costa and Lucas Meyer dos Santos
for developing and debugging the first version of the program.
This research is part of FAPESP ``Projeto Tem\'atico" Grant number
90/3918-5. Partially supported by CNPq.


%
%
%

\parskip 2mm
\parskip 2mm
\noindent
{\bf References}
\vskip0.5truecm
\parskip 2mm
\parindent 0pt

\noindent
[AK] E. Aarts, J. Korst, {\it Simulated Annealing and Boltzmann Machines}, Wiley,
New York, 1989.

[B] F. Barahona, {\it On the computational complexity of Ising spin-glass
models}, J. of Physics A {\bf 15} (1982), 3241.

[FFG] P. A. Ferrari, A. Frigessi and P. Gonzaga de S\'a, 
{\it Fast approximate
maximum a posteriori restoration of multicolor images}, J. R. Statisti. Soc. B
(1995), 485.

[FP] T. A. Ferryman and S. J. Press, {\it A comparison of nine pixel
classification algorithms}, technical report, University of California, (1997).

[FF] L. R. Ford and D. R. Fulkerson, {\it Flows in networks}, 
Princeton: Princeton University Press. (1962).

[FP]  A. Frigessi and M. Piccioni, {\it Parameter estimation for
  two-dimensional Ising fields corrupted by noise}, Stoch. Proc. Appl. {\bf
  34} (1990) 297-311.

[GJ] M. R. Garey and D. S. Johnson, {\it Computers and intractability: a guide
to the theory of NP-completeness}, San Francisco, Freeman (1979).

[G] D. Geman, {\it Random Fields and Inverse Problems in Imaging}, Lecture
Notes in Mathematics {\bf 1470} (1990), Springer Verlag, New York.

[Gi] B.  Gidas {\it Metropolis type Monte Carlo simulation algorithm and
simulated annealing}, Topics in contemporary probability and its applications,
159--232, Probab. Stochastics Ser., CRC, Boca Raton, FL (1995).

[GPS] D. M. Greig, B. T. Porteous, A. M. Seheult, {\it Exact maximum a
  posteriori estimation for binary images}, J. Royal Statistical Society,
Series B {\bf 51} (1989), 271-279.

[J] M. D. Jubb, {\it Ph.D Thesis}, unpublished, University of Bath, UK
(1989).

[M] P. Martin, {\sl
Potts models and related problems in statistical mechanics.}
Series on Advances in Statistical Mechanics, 5.
World Scientific Publishing Co., Inc., Teaneck, NJ, 1991.

[MW] B. M. McCoy and T. T. Wu, {\it The two-dimensional Ising model}, 
Harvard Univ. Press (1973). 

[R] D. Ruelle, {\it Statistical mechanics: Rigorous results}.
 W. A. Benjamin, Inc., New
York-Amsterdam 1969.

[T] A. Trouv\'e
{\it Cycle decompositions and simulated annealing. } 
SIAM J. Control Optim. 34 (1996), no. 3, 966--986.

\vskip 7truemm
\noindent Instituto de Matem\'atica e Estat\'\i stica --- %
Universidade de S\~ao Paulo \hfill\break
Cx.\ Postal 66.281 --- 05315-970 S\~ao Paulo SP --- Brasil \hfill\break
{\tt<pablo@ime.usp.br>}
{\tt<gubi@ime.usp.br>}
{\tt<neves@ime.usp.br>}

\end{document}